\documentclass[12pt]{amsart}
\usepackage{fullpage}
\usepackage{ulem}
\usepackage{amsmath}
\usepackage{amsthm}
\usepackage{amsbsy}
\usepackage{amsrefs}
\usepackage{amssymb}
\usepackage{hyperref}
\usepackage{amsfonts}

\usepackage{enumitem}
\usepackage{amsfonts}
\usepackage[OT2,T1]{fontenc}
\DeclareSymbolFont{cyrletters}{OT2}{wncyr}{m}{n}
\DeclareMathSymbol{\Sha}{\mathalpha}{cyrletters}{"58}
\numberwithin{equation}{section}
\newtheorem{theorem}{Theorem}[section]
\newtheorem{lemma}[theorem]{Lemma}

\newtheorem{proposition}[theorem]{Proposition}
\newtheorem{corollary}[theorem]{Corollary}
\newtheorem{definition}[theorem]{Definition}

\newtheorem{remark}[theorem]{Remark}
\newtheorem*{conjecture*}{Conjecture}
\makeatletter
\def\imod#1{\allowbreak\mkern10mu({\operator@font mod}\,\,#1)}
\makeatother

\newcommand{\Q}{\mathbb{Q}}

    \DeclareFontFamily{U}{wncy}{}
    \DeclareFontShape{U}{wncy}{m}{n}{<->wncyr10}{}
    \DeclareSymbolFont{mcy}{U}{wncy}{m}{n}
    \DeclareMathSymbol{\Sh}{\mathord}{mcy}{"58} 
\title[Bounded Gaps Between Primes in Chebotarev Sets]{Bounded Gaps Between Primes in Chebotarev Sets}
\author{Jesse Thorner}
\date{\today}
\begin{document}
\maketitle
\vspace{-.3in}
\begin{abstract}
A new and exciting breakthrough due to Maynard establishes that there exist infinitely many pairs of distinct primes $p_1,p_2$ with $|p_1-p_2|\leq 600$ as a consequence of the Bombieri-Vinogradov Theorem.  In this paper, we apply his general method to the setting of Chebotarev sets of primes.  We study applications of these bounded gaps with an emphasis on ranks of prime quadratic twists of elliptic curves over $\Q$, congruence properties of the Fourier coefficients of normalized Hecke eigenforms, and representations of primes by binary quadratic forms.
\end{abstract}

\section{Introduction}
The long-standing twin prime conjecture states that there are infinitely many primes $p$ such that $p+2$ is also prime.  The fact that there is a large amount of numerical evidence supporting the twin prime conjecture is fascinating, considering that the Prime Number Theorem tells us that on average, the gap between consecutive primes $p_1,p_2$ is about $\log(p_1)$.  A resolution to the twin prime conjecture seems beyond the reach of current methods.  The next best result for which one could hope is that there are bounded gaps between primes; that is, there exist a constant $C>0$ and infinitely many pairs of distinct primes $p_1,p_2$ satisfying $|p_1-p_2|\leq C$.

In \cite{GPY}, Goldston, Pintz, and Y{\i}ld{\i}r{\i}m developed the ``GPY method", which led to the proof that for any $\epsilon>0$, there exist infinitely many pairs of distinct primes $p_1,p_2$ satisfying
\[
|p_1-p_2|<\epsilon\log(p_1).
\]
This method relies heavily on the distribution of primes in arithmetic progressions.  We say that the primes have {\it level of distribution} $\theta$ if for any fixed $A>0$, we have
\begin{equation}
\label{bv}
\sum_{q\leq x^{\theta}}\max_{y\leq x}\max_{a\in(\mathbb{Z}/q\mathbb{Z})^\times}\left|\pi(y;q,a)-\frac{1}{\varphi(q)}\pi(y)\right|\ll \frac{x}{\log(x)^A},
\end{equation}
where
\begin{equation}
\pi(x)=\#\{\text{$p$ prime: $p\leq x$}\},\quad\pi(x;q,a)=\#\{\text{$p$ prime: $p\leq x,p\equiv a\imod q$}\}.
\end{equation}
As an application of the large sieve, Bombieri and Vinogradov established that (\ref{bv}) holds when $0<\theta<\frac{1}{2}$.  It was conjectured by Elliott and Halberstam \cite{EH} that (\ref{bv}) holds when $0<\theta<1$.  The GPY method produces bounded gaps between primes assuming that $\theta>\frac{1}{2}$.  In \cite{Zhang}, Zhang proved that there are infinitely many pairs of distinct primes $p_1,p_2$ satisfying $|p_1-p_2|\leq 7\times10^7$ by finding a suitable modification for (\ref{bv}) which is valid for $\theta>\frac{1}{2}$.  Zhang's work is inspiring but seems difficult to adapt to other settings.

In \cite{maynard}, Maynard proved that there are infinitely many pairs of distinct primes $p_1,p_2$ satisfying $|p_1-p_2|\leq600$.  (Tao developed the underlying sieve theory independently, but arrived at slightly different conclusions.)  This result follows from a dramatic improvement to the GPY method arising from the use of more general sieve weights.  Once we have this improvement, all that one must know in order to obtain bounded gaps between primes is the distribution of primes within the integers (which is given by the Prime Number Theorem) and the fact that the level of distribution $\theta$ of the primes is positive (which is given by the Bombieri-Vinogradov Theorem).

In this paper, we exploit the flexibility in the methods presented in \cite{maynard} to obtain analogous results on bounded gaps between primes in Chebotarev sets $\mathcal{P}$.  These sets are characterized as follows.  Let $K/\mathbb{Q}$ be a Galois extension of number fields with Galois group $G$ and discriminant $\Delta$.  For a prime $p\nmid\Delta$, there corresponds a certain conjugacy class $C\subset G$ consisting of the set of Frobenius automorphisms attached to the prime ideals of $K$ which lie over $p$.  We denote this conjugacy class by the Artin symbol $[\frac{K/\mathbb{Q}}{p}]$.  We say that a subset $\mathcal{P}$ of the primes is a Chebotarev set, or that $\mathcal{P}$ satisfies a Chebotarev condition, if there exists an extension $K/\mathbb{Q}$ and a union of conjugacy classes $C\subset G$ such that $\mathcal{P}$ is a union of sets of the form $\{\textup{$p$ prime}:p\nmid\Delta,[\frac{K/\mathbb{Q}}{p}]=C\}$.  The Chebotarev Density Theorem asserts that $\mathcal{P}$ has relative density within the primes that is both positive and rational, and a result of Murty and Murty \cite{MM-chebotarev} tells us that we can extend the notion of a positive level of distribution to $\mathcal{P}$ if we omit certain arithmetic progressions.  These two ingredients in conjunction with the sieve developed in \cite{maynard} enable us to prove the existence of bounded gaps between primes in any Chebotarev set.

By the Kronecker-Weber Theorem, if $K/\mathbb{Q}$ is an abelian extension, then $\mathcal{P}$ is determined by congruence conditions.  Thus finding bounded gaps between primes in Chebotarev sets determined by abelian extensions is equivalent to finding bounded gaps between primes in arithmetic progressions, which is proven in \cite{granville} using a combinatorial argument.  In this paper, we handle the nonabelian extensions, proving a complete characterization of bounded gaps between primes in Chebotarev sets.
\begin{theorem}
\label{main-theorem}
Let $K/\mathbb{Q}$ be a Galois extension of number fields with Galois group $G$ and discriminant $\Delta$, and let $C$ be a conjugacy class of $G$.  Let $\mathcal{P}$ be the set of primes $p\nmid\Delta$ for which $[\frac{K/\mathbb{Q}}{p}]=C$.
\begin{enumerate}
\item If $K/\mathbb{Q}$ is a nonabelian extension, then there exist infinitely many pairs of distinct primes $p_1,p_2\in \mathcal{P}$ such that $|p_1-p_2|\leq825(\frac{|G|^2\Delta}{|C|\varphi(\Delta)})^3\exp(\frac{|G|^2\Delta}{|C|\varphi(\Delta)})$.
\item If $K/\mathbb{Q}$ is an abelian extension, let $q$ be the smallest positive integer so that $K\subset \mathbb{Q}(e^{2\pi i/q})$.  There exist infinitely many pairs of distinct primes $p_1,p_2\in\mathcal{P}$ such that $|p_1-p_2|\leq600q$.
\end{enumerate}
\end{theorem}

We use Theorem \ref{main-theorem} to prove several results in algebraic number theory.  The first two of these are immediate.

\begin{corollary}
Let $K/\mathbb{Q}$ be a Galois extension of number fields with ring of integers $\mathcal{O}_K$.  There exist a constant $c(K)>1$ and infinitely many pairs of non conjugate prime ideals $\mathfrak{a},\mathfrak{b}\subset\mathcal{O}_K$ such that $|N_{K/\mathbb{Q}}\mathfrak{a}-N_{K/\mathbb{Q}}\mathfrak{b}|\leq c(K)$.
\end{corollary}

Let $f\in\mathbb{Z}[x]$ be monic polynomial of degree $d$ and discriminant $\Delta$ that is irreducible over $\mathbb{Q}$, and let $G$ be the permutation representation of the Galois group of $f$.  Let $p\nmid\Delta$ be a prime, let $1\leq r\leq d$, and suppose that $f\equiv\prod_{i=1}^r f_i\imod p$ with the $f_i$ distinct irreducible polynomials in $(\mathbb{Z}/p\mathbb{Z})[x]$ of degree $n_i$.  Then $G$ contains a permutation $\sigma_p$ that is a product of disjoint cycles of length $n_i$; we call the cycle type of $\sigma_p$ the {\it factorization type} of $f\bmod p$.

\begin{corollary}
\label{factor}
Assume the above notation.  Let $f\in\mathbb{Z}[x]$ be an irreducible monic polynomial.  There exists a constant $c(f)>1$ and infinitely pairs of distinct primes $p_1,p_2$ such that
\begin{enumerate}
\item $|p_1-p_2|\leq c(f)$.
\item $f\bmod p_1$ and $f\bmod p_2$ have the same factorization type.
\end{enumerate}
\end{corollary}

Theorem \ref{main-theorem} has many interesting applications to the theory of elliptic curves.  Let $E/\mathbb{Q}$ be an elliptic curve, and let $E_d/\mathbb{Q}$ denote the quadratic twist of $E$ by $d$.  We denote the rank of the group of  $\mathbb{Q}$-rational points $E(\mathbb{Q})$ by $\textup{rk}(E)$.  Our applications are related to the following conjecture due to Silverman regarding $\textup{rk}(E_{\pm p})$ when $p$ is prime.
\begin{conjecture*}
\label{silverman}
There are infinitely many primes $p$ for which $\textup{rk}(E_{p})=0$ or $\textup{rk}(E_{-p})=0$, and there are infinitely many primes $\ell$ for which $\textup{rk}(E_{\ell})>0$ or $\textup{rk}(E_{-\ell})>0$.
\end{conjecture*}

In light of Silverman's conjecture, we prove the following result for certain ``good" elliptic curves, which is related to the rank zero component of Silverman's conjecture.

\begin{theorem}
\label{2-trivial-gaps}
Let $E/\mathbb{Q}$ be a ``good" elliptic curve (see Definition \ref{good}).  There exist a constant $c(E)>1$, $\epsilon\in\{-1,1\}$, and infinitely many pairs of distinct primes $p_1,p_2$ such that
\begin{enumerate}
\item $|p_1-p_2|\leq c(E)$.
\item $\textup{rk}(E_{\epsilon p_1})=\textup{rk}(E_{\epsilon p_2})=0.$
\end{enumerate}
\end{theorem}

In light of recent results by Coates, Li, Ye, and Zhai \cite{coates}, we use our results to study ranks of twists of the elliptic curve $E=X_0(49)$, whose minimal Weierstrass equation is given by $E:y^2+xy=x^3-x^2-2x-1.$  Let $p>7$ be a prime such that $p\equiv3\imod4$ and $p$ is inert in the field $\mathbb{Q}(\sqrt{-7})$.  For $k\geq0$, let $q=\prod_{i=1}^k q_i$ be a product of distinct primes $q_i\neq p$, each of which splits completely in $\mathbb{Q}(E[4])$.  Suppose further that the ideal class group of $\mathbb{Q}(\sqrt{-pq})$ has no element of order 4.  Under these hypotheses, Coates, Li, Yian, and Zhai prove that the Hasse-Weil $L$-function $L(E_{-pq},s)$ has a simple zero at $s=1$, $\textup{rk}(E_{-pq})=1$, and the Shafarevich-Tate group $\Sha(E_{-pq}/\mathbb{Q})$ is finite of odd cardinality.  They predict that every elliptic curve should satisfy a property similar to this.  We prove the following.

\begin{theorem}
\label{coates-gaps}
Let $E=X_0(49)$.  There exist infinitely many pairs of distinct primes $p_1,p_2$ such that
\begin{enumerate}
\item $|p_1-p_2|\leq 16800$.
\item Both $L(E_{-p_1},s)$ and $L(E_{-p_2},s)$ have a simple zero at $s=1$, $\textup{rk}(E_{-p_1})=\textup{rk}(E_{-p_2})=1$, and $\Sha(E_{-p_1}/\mathbb{Q})$ and $\Sha(E_{-p_2}/\mathbb{Q})$ are both finite of odd cardinality.
\end{enumerate}
\end{theorem}

A specific elliptic curve for which the entirety of Silverman's conjecture is true is the congruent number elliptic curve $E':y^2=x^3-x$.  We call a positive squarefree integer $d$ {\it congruent} if $d$ is the area of a right triangle with sides of rational length.  It is well-known that $d$ is a congruent number if and only if $E_d'(\mathbb{Q})$ has positive rank.  If $p$ is prime, it is also known \cite{Monsky} that
\begin{itemize}
\item If $p\equiv3\imod 8$, then $\textup{rk}(E_p')=0$.
\item If $p\equiv5\imod 8$, then $\textup{rk}(E_{2p}')=0$.
\item If $p\equiv \text{$5$ or $7$}\imod 8$, then $\textup{rk}(E_p')=1$.
\item If $p\equiv3\imod4$, then $\textup{rk}(E_{2p}')=1$.
\end{itemize}
For such primes, the existence of bounded gaps follows immediately from the second part of Theorem \ref{main-theorem}.  We obtain the following result for twists $E_{p}'(\mathbb{Q})$, but one can easily adapt the statement to suit the twists $E_{2p}'(\mathbb{Q})$.
\begin{theorem}
There exist infinitely many pairs of distinct primes $p_1,p_2$ such that
\begin{enumerate}
\item $|p_1-p_2|\leq 4800$.
\item Either $\textup{rk}(E_{p_1}')=\textup{rk}(E_{p_2})=0$ or $\textup{rk}(E_{p_1}')=\textup{rk}(E_{p_2}')=1$.
\end{enumerate}
In particular, we have bounded gaps between congruent primes and between non-congruent primes.
\end{theorem}

We also have applications to congruence conditions satisfied by the Fourier coefficients of normalized Hecke eigenforms, i.e. newforms, on congruence subgroups of $\textup{SL}_2(\mathbb{Z})$.
\begin{theorem}
\label{fourier-gaps}
Let $f(z)=\sum_{n=1}^\infty a_f(n)q^n\in S_k^{\textup{new}}(\Gamma_0(N),\chi)\cap\mathbb{Z}[[q]]$ be a newform of even weight $k\geq2$, and let $d$ be a positive integer.  There exist a constant $c(d,f)>1$ and infinitely many pairs of distinct primes $p_1,p_2$ such that
\begin{enumerate}
\item $|p_1-p_2|\leq c(d,f)$.
\item $a_f(p_0)\equiv a_f(p_1)\equiv a_f(p_2)\imod d$.
\end{enumerate}
In particular, we have bounded gaps between primes $p$ satisfying $a_f(p)\equiv0\pmod d$.
\end{theorem}
As an application of Theorem \ref{fourier-gaps}, let
\[
f(z)=q\prod_{n=1}^\infty(1-q^n)^{24}=\sum_{n=1}^\infty\tau(n)q^n\in S_{12}^{\textup{new}}(\Gamma_0(1),\chi_{\textup{triv}}),
\]
where $\tau$ is the Ramanujan tau function.  In this case, we have bounded gaps between primes $p$ for which $\tau(p)\equiv0\imod d$ for any positive integer $d$.  If $k=2$, then $f$ is the newform associated to an elliptic curve $E/\mathbb{Q}$ with conductor $N$.  In this case, $a_f(p)=p+1-\#E(\mathbb{F}_p)$, and we have bounded gaps between primes $p$ for which $\#E(\mathbb{F}_p)\equiv p+1\imod d$ for any positive integer $d$.

Finally, we consider primes represented by binary quadratic forms.  Let $Q(x,y)=ax^2+bxy+cy^2\in\mathbb{Z}[x,y]$ be a primitive, positive-definite quadratic form with discriminant $D=b^2-4ac<0$.  It is known that the primes represented by $Q$ form a Chebotarev set.  The proportion of primes that are represented by $Q$ is either $\frac{1}{h(D)}$ or $\frac{1}{2h(D)}$, where $h(D)$ is the class number for quadratic forms $Q$ of discriminant $D$.  If $K=\Q(\sqrt{D})$, $\mathcal{O}$ is the order of the discriminant $D$, and $L$ is the ring class field of $\mathcal{O}$, then the Chebotarev condition satisfied by these primes is in the extension $L/\Q$.  (See \cite{cox} for further discussion.)
\begin{theorem}
\label{binary-gaps}
Let $Q(x,y)=ax^2+bxy+cy^2\in\mathbb{Z}[x,y]$ be a primitive, positive-definite quadratic form with $b^2-4ac<0$.  There exist a constant $c(Q)>1$ and infinitely many pairs of distinct primes $p_1,p_2$ such that
\begin{enumerate}
\item $|p_1-p_2|\leq c(Q)$.
\item Both $p_1$ and $p_2$ are represented by $Q$.
\end{enumerate}
In particular, if $n$ is a positive integer, then there are bounded gaps between primes of the form $x^2+ny^2$.
\end{theorem}

\begin{remark}
\label{k-tuples}
The work of Maynard and Tao also tells us that
\[
\liminf_{n\to\infty}(p_{m+n}-p_n)\ll m^3 e^{4m}.
\]
This extends to Chebotarev sets as well.  If $\{q_1,q_2,q_3,\ldots\}$ is an ordered Chebotarev set of primes corresponding to a nonabelian extension $K/\mathbb{Q}$, then
\[
\liminf_{n\to\infty}(q_{m+n}-q_n)\ll m^3 \exp\left(\frac{m |G|^2\Delta}{|C|\varphi(\Delta)}\right).
\]
The abelian case and all of the above applications have similar results.  The proof is similar to the case of pairs of primes presented here.
\end{remark}

\subsection*{Acknowledgements.}

I would like to thank Tristan Freiberg and Ethan Smith for their useful comments and suggestions on this paper.  I would particularly like to thank James Maynard for spending time discussing his exciting recent work and answering my many questions, as well as Ken Ono for his support and for suggesting this project.  Numerical computations were performed using Mathematica 9.

\section{Notation}

If $f,g$ are real-valued functions, we say that $f(x)=O(g(x))$ if $\limsup_{x\to\infty}|f(x)/g(x)|<\infty$, and $f(x)=o(g(x))$ if $\lim_{x\to\infty}f(x)/g(x)=0$.  Additionally, we use the notation $f(x)\ll g(x)$ to mean the same as $f(x)=O(g(x))$.  We let $\epsilon>0$ be sufficiently small real numbers, where sufficiency will be obvious by context.  The values of $\epsilon$ may vary at each occurrence.  We let $N$ denote a large positive integer, and all asymptotic notation refers to behavior as $N$ tends to infinity.  All sums, products, and suprema are taken over variables in either the primes, which we denote as $\mathbb{P}$, or the positive integers unless otherwise noted.  Any constants implied by the asymptotic notation $o,O,$ or $\ll$ will not depend on $N$, but may depend on $k$, $\mathcal{H}$, $\epsilon$, or other quantities specified in the paper.

We will let $k$ be a fixed positive integer.  The set $\mathcal{H}=\{h_1,\ldots,h_k\}\subset\mathbb{Z}$ will always be admissible; that is, for every prime $p$, the set $\{h_1\bmod p,\ldots,h_k\bmod p\}$ does not contain all residue classes of $\mathbb{Z}/p\mathbb{Z}$.  The functions $\varphi$, $\tau_r(n)$, and $\mu$ refer to the Euler totient function, the number of representations of $n$ as a product of $r$ positive integers, and the M{\"o}bius function, respectively.  We let $p$ be a rational prime; given a set $\mathcal{P}\subset \mathbb{P}$, $p_n$ will denote the $n$-th prime of $\mathcal{P}$.  We let $\#S$ or $|S|$ denote the cardinality of a finite set $S$.  For any $x\in\mathbb{R}$, we let $\lfloor x\rfloor=\max\{a\in\mathbb{Z}:a\leq x\}$ and $\lceil x\rceil=\min\{a\in\mathbb{Z}:a\geq x\}$.  We let $(a,b)=\gcd(a,b)$ and $[a,b]=\text{lcm}(a,b)$.

\section{Bounded Gaps Between Primes}

The variant of the Selberg sieve developed in \cite{maynard} eliminates the $\theta>\frac{1}{2}$ barrier to achieving bounded gaps between primes that the original GPY method encountered.  By studying the proof of the following theorem, it is clear that we obtain bounded gaps between primes as long as $\theta>0$, a condition which is guaranteed by the Bombieri-Vinogradov Theorem.
\begin{theorem}[Maynard]
\label{maynard-tao}
There are infinitely many pairs of distinct primes $p_1,p_2$ with $|p_1-p_2|\leq600.$
\end{theorem}
We provide a brief outline the important components of the proof given in \cite{maynard}.  For a fixed admissible set $\mathcal{H}=\{h_1,\ldots,h_k\}$, we consider the sums
\begin{align}
S_1(N)&=\sum_{\substack{N\leq n< 2N \\ n\equiv v_0\imod W}}w_n,\\
S_2(N)&=\sum_{\substack{N\leq n< 2N \\ n\equiv v_0\imod W}}\sum_{i=1}^k\chi_{\mathbb{P}}(n+h_i)w_n,\\
S(N,\rho)&=S_2(N)-\rho S_1(N),
\end{align}
where $w_n$ are nonnegative weights, $\rho>0$, $\chi_{\mathbb{P}}$ is the indicator function of the primes, and
\begin{equation}
\label{W}
W=\prod_{p\leq D_0}p,\qquad D_0=\log\log\log(N).
\end{equation}
By the Prime Number Theorem, $W\ll \log\log(N)^{2}.$

The goal is to show that $S(N,\rho)>0$ for all sufficiently large $N$.  This would imply that for infinitely many $N$, there exists $n\in[N,2N)$ for which at least $\lfloor \rho+1\rfloor$ of the $n+h_i$ are prime, establishing an infinitude of intervals of bounded length  containing $\lfloor \rho+1\rfloor$ primes.  

Let $\mathcal{R}_k=\{\vec{x}\in[0,1]^k:\sum_{i=1}^k x_i\leq 1\}$, let $F:[0,1]^k\to\mathbb{R}^k$ be a infinitely differentiable function supported on $\mathcal{R}_k$, and let $R=N^{\theta/2-\epsilon}$.  The weights $w_n$ are of the form
\begin{equation}
w_n=\left(\sum_{d_i\mid n+h_i\forall i}\lambda_{d_1,\ldots,d_k}\right)^2,
\end{equation}
where
\begin{equation}
\lambda_{d_1,\ldots,d_k}=\left(\prod_{i=1}^k \mu(d_i)d_i\right)\sum_{\substack{r_1,\ldots,r_k \\ d_i\mid r_i~\forall i \\ (r_i,W)=1~\forall i}}\frac{\mu\left(\prod_{i=1}^k r_i\right)^2}{\prod_{i=1}^k\varphi(r_i)}F\left(\frac{\log(r_1)}{\log(R)},\ldots,\frac{\log(r_k)}{\log(R)}\right).
\end{equation}
Setting $d=\prod_{i=1}^k d_i$, we choose $\lambda_{d_1,\ldots,d_k}$ to be supported when $d<R$, $(d,W)=1$, and $\mu(d)^2=1$.

First, estimates on the sums $S_1(N)$ and $S_2(N)$ are established.
\begin{proposition}
\label{maynard-prop-1}
Let $\mathbb{P}$ have level of distribution $\theta>0$.  Let $F:[0,1]^k\to\mathbb{R}$ be a fixed infinitely differentiable function supported on $\mathcal{R}_k$.  We have
\begin{align*}
S_1(N)&=(1+o(1))\frac{\varphi(W)^k N\log(R)^k}{W^{k+1}}I_k(F),\\
S_2(N)&=(1+o(1))\frac{\varphi(W)^k N\log(R)^{k}}{W^{k+1}}\frac{\log(R)}{\log(N)}\sum_{i=1}^k J_k^{(i)}(F),
\end{align*}
provided that $I_k(F)\neq0$ and $J_k^{(i)}(F)\neq0$ for each $i$, where
\begin{align*}
I_k(F)&=\int_0^1\cdots\int_0^1 F(t_1,\ldots,t_k)^2~dt_1\cdots dt_k,\\
J_k^{(i)}(F)&=\int_0^1\cdots\int_0^1\left(\int_0^1 F(t_1,\ldots,t_k)~dt_i\right)^2~dt_1\cdots dt_{i-1}dt_{i+1}\cdots dt_k.
\end{align*}
\begin{proof}
This is proven in Sections 5 and 6 of \cite{maynard}.
\end{proof}
\end{proposition}
Following the GPY method, we want $S_2(N)-\rho S_1(N)$ to be positive for all sufficiently large $N$, ensuring that for infinitely many $n$, several of the $n+h_i$ are prime.  The following proposition states this formally.
\begin{proposition}
\label{maynard-prop-2}
Let $\mathbb{P}$ have level of distribution $\theta>0$, and let $\mathcal{H}=\{h_1,\ldots,h_k\}$ be an admissible set.  Let $\mathcal{S}_k$ denote the set of infinitely differentiable functions $F:[0,1]^k\to\mathbb{R}$ supported on $\mathcal{R}_k$ with $I_k(F)\neq0$ and $J_k^{(i)}(F)\neq0$ for each $i$.  Define
\[
M_k=\sup_{F\in \mathcal{S}_k}\frac{\sum_{i=1}^k J_k^{(i)}(F)}{I_k(F)},\qquad r_k=\left\lceil \frac{\theta M_k}{2}\right\rceil.
\]
There are infinitely many $n$ such that at least $r_k$ of the $n+h_i$ are prime.  Furthermore, if $p_n$ is the $n$-th prime, then
\[
\liminf_{n\to\infty}(p_{n+r_k-1}-p_{n})\leq\max_{1\leq i<j\leq k}(h_i-h_j).
\]
\end{proposition}
\begin{proof}
This is proven in Section 4 of \cite{maynard}.
\end{proof}

All that remains is to find a suitable lower bound for $M_k$.
\begin{proposition}
\label{maynard-prop-3}
We have $M_{105}>4$.  For sufficiently large $k$,
\[
M_k>\log(k)-2\log(\log(k))-2.
\]
\end{proposition}
\begin{proof}
This is proven in Sections 7 and 8 of \cite{maynard}.
\end{proof}
The exact manner in which these propositions are put together is outlined in Section 4 of \cite{maynard}.  We emulate those arguments in the next section.

\section{Proof of Theorem \ref{main-theorem}}

One fascinating aspect of the proof of Theorem \ref{maynard-tao} is how adaptable it is to exploring bounded gaps between primes in special subsets of the primes.  In this section, we will modify the proof to obtain a version applicable to sets of primes satisfying a Chebotarev condition.

Let $K/\mathbb{Q}$ be a Galois extension of number fields with Galois group $G$ and discriminant $\Delta$, and let $C$ be a conjugacy class of $G$.  Let
\begin{equation}
\mathcal{P}=\left\{\text{$p$ prime}:p\nmid\Delta,\left[\frac{K/\mathbb{Q}}{p}\right]=C\right\},
\end{equation}
where $[\frac{K/\mathbb{Q}}{\cdot}]$ is the Artin symbol, and define
\begin{align}
\pi_{\mathcal{P}}(N)&=\#\{N\leq p< 2N:p\in\mathcal{P}\},\\
\pi_{\mathcal{P}}(N;q,a)&=\#\{N\leq p< 2N:p\in\mathcal{P}:p\equiv a\imod q\}.
\end{align}

We say that $\mathcal{P}$ has level of distribution $\theta$ if there exists a fixed positive integer $M$ such that for any fixed $A>0$,
\begin{equation}
\label{MM-BV}
\sum_{\substack{q\leq N^{\theta} \\ (q,M)=1}}\max_{y\leq N}\max_{a\in(\mathbb{Z}/q\mathbb{Z})^\times}\left|\pi_{\mathcal{P}}(y;q,a)-\frac{1}{\varphi(q)}\pi_{\mathcal{P}}(y)\right|\ll \frac{N}{\log(N)^A}.
\end{equation}

\begin{lemma}
\label{BV}
Assume the above notation.  Let $\delta=|C|/|G|$.  Suppose that $|G|\geq4$.
\begin{enumerate}
\item We have $\pi_{\mathcal{P}}(N)=\delta N/\log(N)+O(N/\log(N)^2).$
\item Equation \ref{MM-BV} holds when $M=\Delta$ and $0<\theta<2/|G|$.
\end{enumerate}
\end{lemma}
\begin{proof}
The first part is the Chebotarev Density Theorem with error term.  The second part follows from the main result in \cite{MM-chebotarev}.
\end{proof}

In order to use the second part of Lemma \ref{BV}, we must modify the work in the previous section.  Let $W$ be defined as in (\ref{W}), and let $\mathcal{H}=\{h_1,\ldots,h_k\}$ be admissible.  For a positive integer $n$, let $\textup{rad}(n)=\prod_{p\mid n}p$. Define
\begin{equation}
\label{UV-def}
\det(\mathcal{H})=\prod_{i\neq j}(h_i-h_j),\qquad U=W/\textup{rad}(\Delta).
\end{equation}
By the Chinese Remainder Theorem and the admissibility of $\mathcal{H}$, there exists an integer $u_0$ satisfying $(\prod_{i=1}^k(u_0+h_i),U)=1$.  Instead of the restriction $n\equiv v_0\imod W$, we use $n\equiv u_0\imod U$.  We note that when $N$ is sufficiently large, $\textup{rad}(\Delta\det(\mathcal{H}))$ divides $W$.  As in the previous section, $\lambda_{d_1,\ldots,d_k}$ will be supported when
\begin{equation}
\label{support}
d=\prod_{i=1}^k d_i<R,\quad (d,W)=1,\quad\mu(d)^2=1,\quad(d_i,d_j)=1\text{ for all $i\neq j$}.
\end{equation}
Therefore, if $N$ is sufficiently large, then (\ref{W}), (\ref{UV-def}), and (\ref{support}) tell us
\begin{align}
\label{support-nonzero}
\lambda_{d_1,\ldots,d_k}\neq0\text{ implies that }\prod_{1\leq i<j\leq k}(d_i,d_j)=(d,U\Delta\det(\mathcal{H}))=1.
\end{align}

Define
\begin{align}
S_1(N,\mathcal{P})&=\sum_{\substack{N\leq n< 2N \\ n\equiv u_0\imod {U}}}\left(\sum_{d_i\mid n+h_i\forall i}\lambda_{d_1,\ldots,d_k}\right)^2,\\
S_2^{(m)}(N,\mathcal{P})&=\sum_{\substack{N\leq n< 2N \\ n\equiv u_0\imod {U}}}\chi_{\mathcal{P}}(n+h_m)\left(\sum_{d_i\mid n+h_i\forall i}\lambda_{d_1,\ldots,d_k}\right)^2,\\
S_2(N,\mathcal{P})&=\sum_{i=1}^k S_2^{(i)}(N,\mathcal{P}),\\
S(N,\rho,\mathcal{P})&=S_2(N,\mathcal{P})-\rho S_1(N,\mathcal{P}),
\end{align}
where $\rho>0$.  For a fixed $\theta>0$ satisfying (\ref{MM-BV}), let $R=N^{\theta/2-\epsilon}$.  We have the following estimate $S_1(N,\mathcal{P})$.
\begin{proposition}
\label{S1}
Assume the above notation.  If $\mathcal{P}$ has level of distribution $\theta>0$, then
\[
S_1(N,\mathcal{P})=(1+o(1))\textup{rad}(\Delta)\frac{\varphi(W)^k N \log(R)^k}{W^{k+1}}I_k(F),
\]
where $I_k(F)$ is defined in Proposition \ref{maynard-prop-1}.
\end{proposition}
\begin{proof}
The only difference between $S_1$ from Proposition \ref{maynard-prop-1} and $S_1(\mathcal{P})$ is that instead of the condition $n\equiv v_0\imod W$, we have $n\equiv u_0\imod U$.  Following the proof of Lemma 5.1 in \cite{maynard}, we will alleviate $S_1(\mathcal{P})$ of any conditions in the sums that depend on $U$.  Then the Selberg sieve manipulations and analysis from \cite{maynard} will give us the desired estimates.

Expanding the square gives us
\[
S_1(N,\mathcal{P})=\sum_{\substack{d_1,\ldots,d_k \\ e_1,\ldots,e_k}}\lambda_{d_1,\ldots,d_k}\lambda_{e_1,\ldots,e_k}\sum_{\substack{N\leq n<2N \\ n\equiv u_0\imod U \\ [d_i,e_i]\mid n+h_i\forall i}}\chi_{\mathcal{P}}(n+h_m).
\]
We now show that we can write the conditions $n\equiv u_0\imod U$ and $[d_i,e_i]\mid n+h_i$ for all $i$ as a single congruence condition.  Let $d=\prod_{a=1}^k d_a$ and $e=\prod_{a=1}^k e_a$.  If $\mu(d)^2=0$ or $\mu(e)^2=0$, then $\lambda_{d_1,\ldots,d_k}\lambda_{e_1,\ldots,e_k}=0$ by (\ref{support}).  Thus we assume that $\mu(d)^2=\mu(e)^2=1$, so each $d_i$ and $e_i$ is squarefree.  With each $d_i,e_i$ squarefree, we consider the following two cases:
\begin{enumerate}
\item If a prime $p$ divides $(U,[d_i,e_i])$ for some $i$, then $p\mid(U,d_i)$ or $p\mid(U,e_i)$.  Thus $p\mid(d,U)$ or $p\mid(e,U)$, and $\lambda_{d_1,\ldots,d_k}\lambda_{e_1,\ldots,e_k}=0$ by (\ref{support-nonzero}).
\item If a prime $p$ divides $([d_i,e_i],[d_j,e_j])$ for some $i\neq j$, then 
\[
\text{$p\mid d~$ or $~p\mid e$,}\qquad p\mid n+h_i,\qquad\text{and}\qquad p\mid n+h_j.
\]
Thus $p\mid(d,h_i-h_j)$ or $p\mid(e,h_i-h_j)$.  Therefore, $p\mid(d,\det(\mathcal{H}))$ or $p\mid(e,\det(\mathcal{H}))$, and $\lambda_{d_1,\ldots,d_k}\lambda_{e_1,\ldots,e_k}=0$ by (\ref{support-nonzero}).
\end{enumerate}
Using the Chinese Remainder Theorem, we conclude that the inner sum can be written as a sum over a single residue class modulo $q=U\prod_{i=1}^k [d_i,e_i]$ when $U$ and each $[d_i,e_i]$ are pairwise coprime, in which case the inner sum is $N/q+O(1)$.  Otherwise, $\lambda_{d_1\ldots,d_k}\lambda_{e_1,\ldots,e_k}=0$.  Using Lemma 5.1 of \cite{maynard} and (\ref{UV-def}), we have
\begin{align*}
S_1(N,\mathcal{P})&=\frac{N}{U}\sideset{}{'}\sum_{\substack{d_1,\ldots,d_k \\ e_1,\ldots,e_k}}\frac{\lambda_{d_1,\ldots,d_k}\lambda_{e_1,\ldots,e_k}}{\prod_{i=1}^k [d_i,e_i]}+O\left(\sideset{}{'}\sum_{\substack{d_1,\ldots,d_k \\ e_1,\ldots,e_k}}|\lambda_{d_1,\ldots,d_k}\lambda_{e_1,\ldots,e_k}|\right)\\
&=\textup{rad}(\Delta)\frac{N}{W}\sideset{}{'}\sum_{\substack{d_1,\ldots,d_k \\ e_1,\ldots,e_k}}\frac{\lambda_{d_1,\ldots,d_k}\lambda_{e_1,\ldots,e_k}}{\prod_{i=1}^k [d_i,e_i]}+O(\lambda_{\max}^2R^2\log(R)^{2k}),
\end{align*}
where $\lambda_{\max}=\sup_{d_1,\ldots,d_k}|\lambda_{d_1,\ldots,d_k}|$ and $\sum'$ denotes the restriction that $U$ and each $[d_i,e_i]$ are pairwise coprime and each $d_i,e_i$ is squarefree.  If a prime $p$ divides $([d_i,e_i],U)$ for some $i$, then we have already shown that $\lambda_{d_1,\ldots,d_k}\lambda_{e_1,\ldots,e_k}=0$.  Therefore, we may take $\sum'$ to denote the condition that $\prod_{i\neq j}([d_i,e_i],[d_j,e_j])=1,$ which is a condition that is independent of the arithmetic progression containing $n$.  Therefore, the condition $\sum'$ is independent of our modulus $U$, as desired.

We now see that $S_1(N,\mathcal{P})$ is a multiple (depending only on $\Delta)$ of $S_1(N)$ in one of the intermediate steps in Lemma 5.1 of \cite{maynard}.  Therefore, the proposition follows from Lemmata 5.1 and 6.2 of \cite{maynard}.
\end{proof}
We will use the reasoning from the above proof to estimate $S_2(N,\mathcal{P})$.
\begin{proposition}
\label{S2}
Assume the above notation.  Let $K/\mathbb{Q}$ be a nonabelian Galois extension of number fields with Galois group $G$ and discriminant $\Delta$, and let $C$ be a conjugacy class of $G$.  Let $\delta=|C|/|G|$.  If the primes in $\mathcal{P}$ have level of distribution $\theta>0$, then
\[
S_2(N,\mathcal{P})=(1+o(1))\delta\varphi(\textup{rad}(\Delta))\frac{\log(R)}{\log(N)}\frac{\varphi(W)^kN\log(R)^{k}}{W^{k+1}} \sum_{i=1}^k J_k^{(i)}(F),
\]
where $J_k^{(i)}(F)$ is defined in Proposition \ref{maynard-prop-1}.
\end{proposition}
\begin{proof}
The desired result follows from estimating each $S_2^{(m)}(N,\mathcal{P})$ for each $1\leq m\leq k$.  Expanding the square gives us
\[
S_2^{(m)}(N,\mathcal{P})=\sum_{\substack{d_1,\ldots,d_k \\ e_1,\ldots,e_k}}\lambda_{d_1,\ldots,d_k}\lambda_{e_1,\ldots,e_k}\sum_{\substack{N\leq n<2N \\ n\equiv u_0\imod U \\ [d_i,e_i]\mid n+h_i\forall i}}\chi_{\mathcal{P}}(n+h_m).
\]
As with $S_1(N,\mathcal{P})$, the inner sum can be written as a sum over a single residue class $a_m$  modulo $q=U\prod_{i=1}^k [d_i,e_i]$ when $U$ and each $[d_i,e_i]$ are pairwise coprime, and $\lambda_{d_1,\ldots,d_k}\lambda_{e_1,\ldots,e_k}=0$ otherwise.

Choose integers $d_1,\ldots,d_k,e_1,\ldots,e_k$ such that $\lambda_{d_1,\ldots,d_k}\lambda_{e_1,\ldots,e_k}\neq0$.  Clearly
\[
a_m\equiv u_0\imod U\quad\text{and}\quad [d_i,e_i]\mid a_m+h_i\text{ for all $i$.}
\]
We can conclude from the support of $\lambda_{d_1,\ldots,d_k}$ and our choices of $u_0$ and $a_m$ that
\[
(u_0+h_m,U)=1\quad\text{and}\quad(h_m-h_i,[d_i,e_i])=1\text{ for all $i\neq m$,}
\]
so
\[
(q/[d_m,e_m],a_m+h_m)=1\quad\text{and}\quad [d_m,e_m]\mid a_m+h_m.
\]
Therefore, $(q,a_m+h_m)=1$ if and only if $d_m=e_m=1$.  In this case, the inner sum will have size $\pi_{\mathcal{P}}(N)/\varphi(q)+O(E(N,q))$, where
\[
E(N,q)=\max_{(a,q)=1}\left|\pi_{\mathcal{P}}(N;q,a)-\frac{1}{\varphi(q)}\pi_{\mathcal{P}}(N)\right|.
\]

If $(q,a_m+h_m)\neq1$, then the inner sum equals either 0 or 1.  The inner sum equals 1 if and only if there exists a prime $p$ satisfying $n+h_m=p$ for some $n\in[N,2N)$ with $p\mid q$.  Since $N$ is large, we have $N-|h_m|>\sqrt{N}>R$.  Thus $n+h_m=p$ for some $n\in[N,2N)$ implies that $p>R$, so if $p\mid q$, then $\lambda_{d_1,\ldots,d_k}\lambda_{e_1,\ldots,e_k}=0$.  Thus the inner sum only contributes to $S_2^{(m)}(\mathcal{P})$ when $(q,a_m+h_m)=1$.  We conclude that
\[
S_2^{(m)}(N,\mathcal{P})=\frac{\pi_{\mathcal{P}}(N)}{\varphi(U)}\sideset{}{'}\sum_{\substack{d_1,\ldots,d_k \\ e_1,\ldots,e_k \\ d_m=e_m=1}}\frac{\lambda_{d_1,\ldots,d_k}\lambda_{e_1,\ldots,e_k}}{\prod_{i=1}^k \varphi([d_i,e_i])}+O\left(\sideset{}{'}\sum_{\substack{d_1,\ldots,d_k \\ e_1,\ldots,e_k}}|\lambda_{d_1,\ldots,d_k}\lambda_{e_1,\ldots,e_k}|\cdot E(N,q)\right),
\]
where $q=U\prod_{i=1}^k[d_i,e_i]$ and $\sum'$ denotes the restriction that $U$ and each $[d_i,e_i]$ be pairwise coprime.

We first analyze the error term.  From the support of $\lambda_{d_1,\ldots,d_k}$, we only need to consider squarefree $q<R^2 U\leq N^{\theta-\epsilon}$ satisfying $(q,\Delta)=1$, where $\epsilon>0$ is sufficiently small.  Given a squarefree integer $r$, there are at most $\tau_{3k}(r)$ choices of $d_1,\ldots,d_k,e_1,\ldots,e_k$ for which $r=U\prod_{i=1}^k [d_i,e_i]$.  From Lemma 5.2 of \cite{maynard}, the error term is now
\[
\ll \lambda_{\max}^2\sum_{\substack{r<N^{\theta-\epsilon} \\ (r,\Delta)=1}}\mu(r)^2\tau_{3k}(r)E(N,r),
\]
Using the Cauchy-Schwarz inequality and the trivial bound $E(N,q)\ll N/\varphi(q)$,  the error term is
\[
\ll \lambda_{\max}^2\left(\sum_{\substack{r<N^{\theta-\epsilon} \\ (r,\Delta)=1}}\mu(r)^2\tau_{3k}(r)^2\frac{N}{\varphi(r)}\right)^{1/2}\left(\sum_{\substack{r< N^{\theta-\epsilon} \\ (r,\Delta)=1}} \mu(r)^2 E(N,r)\right)^{1/2}.
\]
It follows from elementary bounds on $\tau_{3k}(r)$ and Lemma \ref{BV} that the error is $\ll \lambda_{\max}^2 N/\log(N)^A$ for any fixed $A>0$, which is also true of the error term in $S_2^{(m)}(N)$ in Lemma 5.2 of \cite{maynard}.

Using (\ref{UV-def}), for any fixed $A>0$, we have
\[
S_2^{(m)}(N,\mathcal{P})=\varphi(\textup{rad}(\Delta))\frac{\pi_{\mathcal{P}}(N)}{\varphi(W)}\sideset{}{'}\sum_{\substack{d_1,\ldots,d_k \\ e_1,\ldots,e_k \\ d_m=e_m=1}}\frac{\lambda_{d_1,\ldots,d_k}\lambda_{e_1,\ldots,e_k}}{\prod_{i=1}^k \varphi([d_i,e_i])}+O(\lambda_{\max}^2 N/\log(N)^A),
\]
where $\sum'$ denotes the restriction that $U$ and each $[d_i,e_i]$ be pairwise coprime.  As in the proof of Proposition \ref{S1}, we can take $\sum'$ to denote the restriction that $\prod_{i\neq j}([d_i,e_i],[d_j,e_j])=1$.  We now see that up to the choice of prime counting function (which results in the factor of $\delta$ in the statement of the proposition), $S_2^{(m)}(N,\mathcal{P})$ is a multiple (depending only on $\delta$ and $\Delta)$ of $S_2^{(m)}(N)$ in one of the intermediate steps in Lemma 5.2 of \cite{maynard}.  Therefore, the proposition follows from Lemmata 5.2 and 6.3 of \cite{maynard} and Lemma \ref{BV}.
\end{proof}

We now modify Proposition \ref{maynard-prop-2} accordingly.
\begin{proposition}
\label{maynard-prop-2-lambda}
Let $\mathcal{H}=\{h_1,\ldots,h_k\}$ be an admissible set, let $\mathcal{P}$ have level of distribution $\theta>0$, and let
\[
M_k=\sup_{F\in \mathcal{S}_K}\frac{\sum_{i=1}^k J_k^{(i)}(F)}{I_k(F)},\qquad r_k=\left\lceil\frac{\delta\theta \varphi(\Delta) M_k}{2\Delta}\right\rceil.
\]
Then there are infinitely many $n$ such that at least $r_k$ of the $n+h_i$ are in $\mathcal{P}$.  Furthermore, if $p_n$ is the $n$-th prime in $\mathcal{P}$, then
\[
\liminf_{n\to\infty}(p_{n+r_k-1}-p_n)\leq \max_{1\leq i<j\leq k}(h_i-h_j).
\]
\end{proposition}
\begin{proof}
We want to show that $S(N,\rho,\mathcal{P})>0$ for all sufficiently large $N$.  Recall that $R=N^{\theta/2-\epsilon}$ for some small $\epsilon>0$.  By the definition of $M_k$, we can choose $F_0\in\mathcal{S}_k$ such that
\[
\sum_{i=1}^k J_k^{(i)}(F_0)>(M_k-\epsilon)I_k(F_0).
\]
Using Propositions \ref{S1} and \ref{S2} and the identity $\frac{\varphi(\textup{rad}(\Delta))}{\textup{rad}(\Delta)}=\frac{\varphi(\Delta)}{\Delta}$, we have
\begin{align*}
&S(N,\rho,\mathcal{P})\\
&=\frac{\varphi(W)^k N\log(R)^k}{W^{k+1}}\left(\frac{\log(R)}{\log(N)}\delta\varphi(\textup{rad}(\Delta))\sum_{i=1}^k J_k^{(i)}(F_0)-\rho \textup{rad}(\Delta) I_k(F_0)+o(1)\right)\\
&\geq\frac{\varphi(W)^kN\log(R)^k I_k(F_0)}{W^{k+1}}\left(\frac{\delta\varphi(\Delta)}{\Delta}\left(\frac{\theta}{2}-\delta\right)(M_k-2\delta)-\rho+o(1)\right).
\end{align*}
Let $\rho=M_k(\frac{\delta\theta\varphi(\Delta)}{2\Delta}-\epsilon).$  By choosing $\delta$ suitably small (depending on $\epsilon$), we have $S(N,\rho,\mathcal{P})>0$ for all sufficiently large $N$.  Thus there are infinitely many $n$ for which at least $\lfloor \rho+1\rfloor$ of the $n+h_i$ are in $\mathcal{P}$.  If $\epsilon$ is sufficiently small, then $\lfloor \rho+1\rfloor=\left\lceil\frac{\delta\theta\varphi(\Delta) M_k}{2\Delta}\right\rceil,$ and we obtain the claimed result.
\end{proof}
It remains to find a suitable lower bound for $M_k$.  Proposition \ref{maynard-prop-3} gives us a lower bound on $M_k$ when $k$ is sufficiently large.  We will now establish the full range of $k$ for which this lower bound holds.   Although this lower bound is far from optimal for $k$ low in the range, the following suffices for the purposes of this paper since $k$ will typically be very large.
\begin{proposition}
\label{maynard-prop-3-lambda}
Let $k\geq213$ be a positive integer.  We have
\[
M_k>\log(k)-2\log(\log(k))-2.
\]
\end{proposition}
\begin{proof}
By the analysis in Section 8 of \cite{maynard}, for some positive constant $A$, we have
\[
M_k\geq A\left(1-\frac{Ae^{A}}{k(1-\frac{A}{e^A-1}-\frac{e^A}{k})^2}\right),
\]
provided that the right hand side is positive.  Let $A=\log(k)-2\log(\log(k))$ as in \cite{maynard}.  With this choice of $A$, the left hand side of the above inequality is bounded below by $\log(k)-2\log(\log(k))-2$ for $k\geq16$.  Since $\log(k)-2\log(\log(k))-2>0$ when $k\geq213$, we have the desired result.
\end{proof}
\begin{proof}[Proof of Theorem \ref{main-theorem}]
The second part is proven in \cite{granville}; it remains to prove the first part.  Suppose that $K/\mathbb{Q}$ is nonabelian.  Since $|G|\geq6$, Lemma \ref{BV} tells us that $\mathcal{P}$ has level of distribution $\theta=2/|G|-\epsilon$.  By Proposition \ref{maynard-prop-3-lambda}, if $k\geq213$ is an integer, then
\begin{equation}
\label{Mk-bound}
\frac{\delta\theta\varphi(\Delta) M_k}{2\Delta}\geq\frac{\delta(2/|G|-\epsilon)\varphi(\Delta)}{2\Delta}(\log(k)-2\log(\log(k))-2).
\end{equation}
Choosing $\epsilon=\frac{2}{k|G|}$ and recalling that $\delta=|C|/|G|$, it follows from a numerical calculation that (\ref{Mk-bound}) is greater than $1$ if
\[
k=125\left\lceil \left(\frac{|G|^2\Delta}{|C|\varphi(\Delta)}\right)^2\exp\left(\frac{|G|^2\Delta}{|C|\varphi(\Delta)}\right)\right\rceil.
\]
Our expression for $k$ is minimized when $\Delta=1$ and $|C|=|G|=6$, in which case $k\geq213$ as required by Proposition \ref{maynard-prop-3-lambda}.  Thus for any admissible set $\mathcal{H}=\{h_1,\ldots,h_k\}$ with $k$ as above, at least $2$ of the $n+h_i$ are in $\mathcal{P}$ for infinitely many integers $n$.  We can choose $h_j=q_{\pi(k)+j}$, where $1\leq j\leq k$ and $q_j$ is the $j$-th prime in $\mathbb{P}$.  For $n\geq6$, we have \cite{Dusart1}
\[
n\log(n)+n\log(\log(n))-n<q_n<n\log(n)+n\log(\log(n)).
\]
Furthermore, if $n\geq355991$, we have \cite{Dusart2}
\[
\frac{n}{\log(n)}\left(1+\frac{1}{\log(n)}\right)\leq\pi(n)\leq\frac{n}{\log(n)}\left(1+\frac{1}{\log(n)}+\frac{2.51}{\log(n)^2}\right).
\]
By a numerical calculation, $q_{\pi(k)+k}-q_{\pi(k)+1}\leq1.6 k\log(k)$ for all $k\geq213$.  Therefore, if $p_n$ is the $n$-th prime of $\mathcal{P}$, then
\begin{align*}
\liminf_{n\to\infty}(p_{n+1}-p_n)&\leq \max_{1\leq i<j\leq k}(h_j-h_i)\\
&\leq 1.6 k\log(k)\\
&\leq825\left(\frac{|G|^2\Delta}{|C| \varphi(\Delta)}\right)^3\exp\left(\frac{|G|^2\Delta}{|C|\varphi(\Delta)}\right).
\end{align*}
\end{proof}

\section{Proofs of Theorems \ref{2-trivial-gaps}, \ref{coates-gaps}, and \ref{fourier-gaps}}

To prove Theorems \ref{2-trivial-gaps}, \ref{coates-gaps}, and \ref{fourier-gaps}, it suffices to prove that the set of primes in each theorem is a Chebotarev set.  The claimed bounds will follow from Theorem \ref{main-theorem}.

\subsection{Ranks of elliptic curves}

To prove Theorem \ref{2-trivial-gaps}, let $E/\mathbb{Q}$ be an elliptic curve with Weierstrass form
\[
E:y^2=x^3+ax^2+bx+c,\quad a,b,c\in\mathbb{Z},
\]
where the discriminant of the cubic is nonzero.  We will assume that $E$ and its points are $\mathbb{Q}$-rational.  If $d$ is a squarefree integer, we define $E_d$ to be the $d$-quadratic twist of $E$ given by
\[
E_d:dy^2=x^3+ax^2+bx+c.
\]
\begin{definition}
\label{good}
Let $E/\mathbb{Q}$ be an elliptic curve without $\mathbb{Q}$-rational 2-torsion.  Following \cite{boxer-diao}, we call $E$ good if $E$ satisfies the following criteria:
\begin{enumerate}
\item The $2$-Selmer rank of $E$ is zero.
\item The discriminant $\Delta$ of $E$ is negative.
\item If $p$ is any prime for which $E$ has bad reduction, then $E$ has multiplicative reduction at $p$, and $v_p(\Delta)$ is odd.
\item $E$ has good reduction at $2$ and the reduction of $E$ modulo 2 has $j$-invariant zero.
\end{enumerate}
\end{definition}
A prototypical example of a good elliptic curve is $E=X_0(11)$, which has Weierstrass form $E:y^2=x^3-4x^2-160x-1264$.

We define a squarefree integer $d$ to be {\it $2$-trivial} for $E$ if $E$ has no rational $2$-torsion modulo $p$ for every odd prime $p\mid d$.  For good elliptic curves, the following is proven in \cite{boxer-diao}.
\begin{theorem}
\label{boxer-diao}
Let $E/\mathbb{Q}$ be a good elliptic curve.  If $d$ is a squarefree $2$-trivial integer for $E$ with $(d,\Delta)=1$, then
\[
\textup{dim}_{\mathbb{F}_2}(\textup{Sel}_2(E_d(\mathbb{Q}))=\begin{cases}
0&\mbox{if $d$ is odd,}\\
1&\mbox{if $d$ is even.}
\end{cases}
\]
In particular, for such odd $d$, we have $\textup{rk}(E_d)=0$.
\end{theorem}

We now prove Theorem \ref{2-trivial-gaps}.

\begin{proof}[Proof of Theorem \ref{2-trivial-gaps}]
We write $E$ in Weierstrass form $E:y^2=f(x)$, where $f(x)=x^3+ax^2+bx+c\in\mathbb{Z}[x]$ has Galois group $G$ and discriminant $\Delta$.  Since $E$ is good, $f$ is irreducible over $\mathbb{Z}$ and $G\cong S_3$.  By the above discussion, the primes $p$ satisfying the hypotheses of Theorem \ref{boxer-diao} are exactly the primes $p\nmid\Delta$ such that $f\bmod p$ is irreducible, that is, the factorization type of $f\bmod p$ corresponds to 3-cycles in $S_3$.  The desired result now follows from Corollary \ref{factor}.
\end{proof}

We will use Theorem 91 of \cite{coates}, which we now state, to prove Theorem \ref{coates-gaps}.

\begin{theorem}
\label{coates-BSD}
Let $E=X_0(49)$.  For $k\geq0$, let $p,q_1,\ldots,q_k$ be prime, and let $N=pq_1q_2\cdots q_k$ be a product of distinct primes satisfying
\begin{enumerate}
\item $p\equiv3\imod4$, $p\neq7$, and $p$ is a quadratic non-residue modulo $7$.
\item $q_1,\ldots,q_k$ split completely in $\mathbb{Q}(E[4])$.
\item The ideal class group $\mathcal{H}_N$ of the field $\mathbb{Q}(\sqrt{-N})$ has no element of order $4$.
\end{enumerate}
Then the Hasse-Weil $L$-function $L(E_{-N},s)$ has a simple zero at $s=1$, $E_{-N}(\mathbb{Q})$ has rank $1$, and the Shafarevich-Tate group of $E_{-N}$ is finite of odd order.
\end{theorem}

\begin{proof}[Proof of Theorem \ref{coates-gaps}]
We consider the case of Theorem \ref{coates-BSD} where $k=0$.  Using the theory of quadratic forms, Gauss proved that if $p\equiv3\imod4$, then $|\mathcal{H}_p|$ is odd.  Thus Theorem \ref{coates-BSD} holds when $N$ is a prime such that  $N\neq7$ such that $N\equiv3\imod 4$ and $N$ is a quadratic non-residue modulo 7.   Every prime $p$ congruent to 3, 19, or 27 modulo 28 satisfies this condition, and the desired result follows from the second part of Theorem \ref{main-theorem}.
\end{proof}

\subsection{Coefficients of newforms}

Following Murty and Murty \cite{MM-fourier}, let $q=e^{2\pi i z}$, and let
\[
f(z)=\sum_{n=1}^\infty a_f(n)q^n\in S_{k}^{\text{new}}(\Gamma_0(N),\chi)\cap\mathbb{Z}[[q]]
\]
be a newform of even weight $k\geq2$ and character $\chi$.  (This forces $\chi$ to be real, and $\chi$ is nontrivial if and only if $f$ has complex multiplication.)  Let $G=\textup{Gal}(\bar{\mathbb{Q}}/\mathbb{Q})$, and let $d$ be a positive integer.  By the work of Deligne, there exists a representation
\[
\rho_d:G\to\textup{GL}_2\left(\prod_{\substack{\textup{$\ell$ prime} \\ \ell\mid d}}\mathbb{Z}/\ell\mathbb{Z}\right)
\]
with the property that if $p\nmid dN$ is prime and $\sigma_p$ is a Frobenius element at $p$ in $G$, then $\rho_d$ is unramified at $p$ and
\[
\textup{tr}\rho_d(\sigma_p)=a_f(p),\quad\det\rho_d(\sigma_p)=\chi(p)p^{k-1}.
\]
Let $\tilde{\rho}_d:G\to \textup{GL}_2(\mathbb{Z}/d\mathbb{Z})$ be the reduction modulo $d$ of $\rho_d$.  Let $H_d$ be the kernel of $\tilde{\rho}_d$, let $K_d$ be the subfield of $\bar{\mathbb{Q}}$ fixed by $H_d$, and let $G_d=\textup{Gal}(K_d/\mathbb{Q})$.  If $q\nmid dN$ is prime, then the condition $a_f(q)\equiv0\pmod d$ means that for any Frobenius element $\sigma_q$ of $q$, $\tilde{\rho}_d(\sigma_q)\in C_d$.  Since $C_d$ contains the image of complex conjugations, $C_d$ is nonempty.

\begin{proof}[Proof of Theorem \ref{fourier-gaps}]
By the preceding discussion, the set of primes $p$ for which $a_f(p)\equiv0\pmod d$ is a Chebotarev set.  By a slight variation of the preceding discussion, we find that for any fixed prime $p_0\nmid dN$, the set of primes $p$ for which $a_f(p)\equiv a_f(p_0)\imod d$ is also a Chebotarev set.  The desired result now follows from Theorem \ref{main-theorem}.
\end{proof}

\bibliography{Paper}
\end{document}